\documentclass[12pt]{amsart}

\usepackage{amsmath,amssymb, amscd, amsfonts}

\def\CC{{\mathbb C}}

\def\QQ{{\mathbb Q}}
\def\DD{{\mathbb D}}
\def\NN{{\mathbb N}}
\def\PP{{\mathbb P}}

\def\bd{{\partial}}

\newcommand {\ga} {\gamma}
\newcommand {\ka} {\kappa}
\newcommand {\kp} {\kappa_+^{\prime}}

\newcommand {\Om} {\Omega}

\newcommand {\xd} {{X^\partial}}
\newcommand {\fd} {{f^\partial}}
\newcommand {\gd} {{g^\partial}}
\newcommand {\hd} {{h^\partial}}

\newcommand {\BOX} {\rule{2mm}{2mm}}

\newcommand {\m} {{\setminus}}
\newcommand {\dv} {{\rm Div}}

\newcommand {\be} {\begin{equation}}
\newcommand {\eeqn} {\end{equation}}
\newcommand {\bea} { \begin{eqnarray}}
\newcommand {\eea} {\end{eqnarray}}
\newcommand {\beas} { \begin{eqnarray*}}
\newcommand {\eeas} {\end{eqnarray*}}
\newcommand {\ra} {\rightarrow}
\newcommand {\dra} {\dashrightarrow}
\newcommand {\hra} {\hookrightarrow}

\newtheorem {lemma} {LEMMA} [section]
\newtheorem {theorem}[lemma]{THEOREM}
\newtheorem {prop}[lemma]{PROPOSITION}
\newtheorem {lem}[lemma]{LEMMA}
\newtheorem {cor}[lemma]{COROLLARY}
\newtheorem {defn}[lemma]{DEFINITION}

\newtheorem {definition}[lemma]{DEFINITION}

\numberwithin{equation}{section}

\begin{document}


\title{Quasiprojective varieties admitting Zariski dense entire
holomorphic curves}
\author{Steven S. Y. Lu
\ \ \ \ Jorg Winkelmann
}

\begin{abstract} Let $X$ be a complex quasiprojective variety.
A result of Noguchi-Winkelmann-Yamanoi shows that 
if $X$ admits a Zariski dense entire 
curve, then its quasi-Albanese map is a fiber space. 
We show that the orbifold structure induced by a properly
birationally equivalent map on the base is special in this case. 
As a consequence, if $X$ is of log-general type with  $\bar q(X)\geq\dim X$, 
then any entire curve is 
contained in a proper subvariety in $X$.
\end{abstract}

\address{%
Steven S. Y. Lu\\
D\'epartement de Math\'ematiques\\
Universit\'e du Qu\'ebec \`a Montr\'eal\\
C.P.~8888, Succursale Centre-ville\\
Montr\'eal, Qc H3C 3P8\\
Canada}
\email{lu.steven@uqam.ca\\}

\address{%
J\"org Winkelmann \\
Mathematisches Institut\\
Universit\"at Bayreuth\\
Universit\"atsstra\ss e 30\\
D-95447 Bayreuth\\
Germany\\
}
\email{jwinkel@member.ams.org}
\maketitle
\section{Introduction and the statement of the main result}

This paper deals with a question of Campana concerning 
the characterization of complex algebraic varieties
that admits Zariski dense entire holomorphic curves. 
This problem for an algebraic surface not of 
log-general type nor a very general algebraic K3 surface
was solved completely by \cite{BL1, BL2}. Campana
in \cite{Ca02} introduced the notion of special varieties
which is a practical way to extend the so far conjectural 
characterization to higher dimensions. Using a recent 
result of Noguchi-Winkelmann-Yamaoi \cite{NWY2}, 
we verify one direction of this characterization here for 
all algebraic varieties whose
quasi-Albanese map is generically finite.

Given a complex projective manifold $\bar X$ with a normal crossing
divisor on it,
we call the pair $X=(\bar X, D)$ a log-manifold. Recall that there is a
locally
free subsheaf of the holomorphic tangent sheaf of $\bar X$ , called the
log-tangent sheaf of $X$, which we denote by $\bar T_X$. It
is the sheaf of holomorphic vector fields leaving $D$ invariant.
Its dual $\bar \Om_X=\bar T^{\vee}_X$ is called the log-cotangent sheaf of
$X$ and $\bar K_X=\det \bar \Om_X$  the log-canonical sheaf of $X$. 
Their sections
are called logarithmic 1-forms, respectively logarithmic volume forms.
Here, and later we will consistently abuse notation
and identify holomorphic vector bundles with their sheaves of sections.
We will abuse the notation further at times and identify a line bundle with
a divisor it corresponds to, for example in the identification
 $\bar K_X=K_{\bar X}(D)= K_{\bar X}+D.$
We first give some proper birational invariants of $\bar X\setminus D$,
which we will also identify with $X$ by a standard abuse of notation.

\begin{definition} With this setup, we define the log-irregularity of $X$ by
$$\bar q(X)=\dim H^0(\bar\Om_X)$$
and we define the log-Kodaira dimension of $X$
by $\bar\ka (X)= \kappa(\bar K_X)$,
where the Kodaira dimension for an invertible sheaf $L$ is given by
$$\ka (L)=\limsup_{m\ra \infty} \frac{ \log \dim H^0(L^{\otimes m})}{\log
m}.$$
We also define, see \cite{Ca02, Lu02},
the essential or the core dimension of $X$ by
$$\kp(X)=\max\{\ p\ |\ L\hra \bar\Omega_X^p\
\text{ is an invertible subsheaf with}\
\ka(L)=p,\ 0\le p\le \dim X\ \}$$

\end{definition}

\noindent
It is an easy fact that the Kodaira dimension of an invertible sheaf
$L$ is invariant under positive tensor powers of $L$ and so
the Kodaira dimension $\ka$ 
makes sense for $\QQ$-invertible sheaves of the form $L(A)$
where $A$ is a $\QQ$-divisor.
We recall also the fact that $\kappa(L)\in \{\infty,0,\dots, \dim(X)\}$
for a $\QQ$-invertible sheaf $L$ and that sections
of powers of $L$, if exist,  define a dominant rational map
$I_L$ to a projective variety of dimension
$\ka(L)$, called the Iitaka fibration of $L$.
We usually allow $I_L$ to be defined on any smooth birational model of
$\bar X$ and choose a model on which $I_L$ is a morphism.
With such a choice, recall that the general (in fact generic)
fibers $F$ of $I_L$  are connected and 
$\ka(I|_F)=0$, see e.g. \cite{Ue}. \\

Let $X_0$ be a quasiprojective variety and $\bar X_0$ a projectivization.
We recall that a log-resolution of $X_0$ is a birational morphism
$r: \bar X\ra \bar X_0$ where
$(\bar X, D)$ is a log-manifold with $D=r^{-1}(\bar X_0\m X_0)$. Such a
resolution
exists by the resolution of singularity theorem. The Hartog extension theorem
allows us to define $\bar q$, $\bar\ka$ and $\kp$ for $X_0$ by taking them
to be
those of a log-resolution.
These are thus proper birational invariants of $X_0$. Here, a proper
birational
map between two quasi-projective varieties are just compositions of proper
birational morphisms and their inverses. Another proper
birational invariant is
given by the (quasi-)Albanese map of $X_0$, which is an algebraic morphism
$$\alpha_{X} : X \ra {\rm Alb}(X)=\colon {\rm Alb}(X_0)$$
defined for any log-resolution $(\bar X,D)$ by line integrals of
logarithmic 1-forms on
$\bar X\setminus D$
with a choice of base point outside $D$ where Alb$(X)$ is a $\bar q(X_0)$
dimensional
semi-abelian variety called the Albanese variety of $X$,  see for example
\cite{NWY1}.
Implicit here is the invariance of ${\rm Alb}(X)$ and the compatibility of
the Albanese map
among log-resolutions. We recall that given a compactification of Alb$(X_0)$,
there exist a
log-resolution of $X_0$ (a compactification of $X_0$ by normal crossing
divisors
in the case $X_0$ is smooth) over which $\alpha_X$ extends to a morphism.
We recall also that a semi-abelian variety is a complex abelian group
$T$ that admits a semidirect product structure via a holomorphic
exact sequence of groups
\begin{equation}\label{2}
0 \ra (\CC^*)^{k} \ra T \stackrel{\pi}{\ra} A \ra 0\ ,
\end{equation}
where $A$ is an abelian variety and $k\ge 0$. It follows that the
algebro-geometric
image of the Albanese map (or the Albanese image) of $X_0$, its dimension
as well as the Albanese variety are proper
birational invariants of $X_0$.

\begin{definition}
We say that $X_0$ is {\bf\em special} if $\kp(X_0)=0$ and that it is
of {\bf\em general type} (or if more precision is required,
of {\bf\em log-general type}) if ${\bar\kappa}(X_0)=\dim X_0$.
\end{definition}

Recall that a holomorphic image of a curve in $X_0$ is called
algebraically degenerate if it is not Zariski dense.
Our main theorems in this paper
is as follows.

\begin{theorem}
Let $X$ be a complex quasi-projective variety with $\bar q(X)\geq \dim X$.
Then every entire holomorphic curve in $X$ is
algebraically degenerate if $X$ is not special. 
Alternatively $X$ admits a Zariski dense 
entire holomorphic curve only if $X$ is special. 
\end{theorem}

\begin{cor} With the same hypothesis on $X$, let $f: \CC\ra X$
be holomorphic and nontrivial. If $X$ is of general type,
then there is a proper subvariety of $X$ containing $f(\CC)$.
\end{cor}

We note that \cite{NWY3} has proved the
same theorem with $\kp$ replaced by $\bar \ka$ but
with the additional hypothesis that the Albanese map of $X$ is proper
and generically finite.
However, without the properness condition for
the Albanese map, the birational condition of $\bar\kappa (X)=0$
is not implied by the
condition that $X$ admits a Zariski dense holomorphic image of $\CC$,
see \cite{DL}.

An important part of this paper is an adaptation 
to the context of special varieties of the results of
Noguchi-Winkelmann-Yamanoi
\cite{NWY1, NWY2,NWY3}
concerning varieties that admit finite maps to semi-Abelian
varieties. All the results on special varieties used for the main 
theorem here are worked out
here from scratch independently of previous sources.
The second author has
spoken about the result on surfaces at a workshop 
at the Fields Institute in 2008 that claimed the connection 
with the characterization by special varieties. This
connection, at least in one direction of the characterization, 
is fully worked out here for
all dimensions. The first author would like to  
thank Gerd Dethloff for valuable discussions on Nevanlinna theory
and especially for the last part of the proof in
proposition~\ref{last} of the paper. He would also like
to thank Fr\'ed\'eric Campana for agreeing on certain
new terminologies introduced in this paper, especially
the use of ``base-special" to characterize a notion introduced
and the accompanying use of ``base-general(-typical)" as an 
alternative for an old notion.

\section{Preliminaries on special varieties}

Throughout this section, let $X$ be
a complex projective manifold and Div$'(X)$ the
set of codimension-one subvarieties of $X$. An orbifold structure on $X$ is a
$\QQ$-divisor of the form $$A=\sum_i (1-1/m_i) D_i$$ where
$1\le m_i\in \QQ\cup\{\infty\}$ and $D_i\in \dv'(X)$ for all $i$.
We denote $X$ with its orbifold structure by $X\m A$ and we set
$K_{X\m A}:= K_X(A)$ to be the orbifold canonical $\QQ$-bundle. We set
$m(A\cap D_i)=m(D_i\cap A)=m_i$ and
call it the multiplicity of the orbifold $X\m A$
(or simply, the orbifold multiplicity) at $D_i$.
We note that the coefficient of $D_i$ in $A$ satisfies
$$0\le1-\frac{1}{m(A\cap D_i)}\leq 1$$
and it vanishes, respectively equals one,
precisely when the corresponding orbifold multiplicity is one,
respectively equals $\infty$. Note  that when $A$ is a (reduced) normal
crossing divisor,
the orbifold $X\m A$ is nothing but a log-manifold $(X,A)$ whose
birational geometry is dictated precisely by the proper birational geometry
of the complement of $A$ in $X$. More generally, when $A_{red}$ is normal
crossing, one can make good geometric sense of the orbifold $X\m A$ via
the usual
branched covering trick (see \cite{Lu02}, see also \cite{Ca08} for a
variant approach)
and we will call such an orbifold smooth.\\

We now define the Kodaira dimension of a rational
map from an orbifold following \cite{Lu02}, c.f. also \cite{Ca02}.
Let $f: X \dra Y$ be a rational map between complex projective manifold
and let $w$ be a rational section of $K_Y$.
If $f$ is dominant, then $f^*w$ defines a rational section of $\Omega_X^m$
 with $m=\dim Y$ and
hence determines in the standard way a saturated rank-one subsheaf $L$ of
$\Omega_X^m$ which is easily seen to be unique in the birational equivalence
class of $f$ (it is even unchanged after composing with a dominant map
from $Y$ to a variety of the same dimension as $Y$).
We recall that a saturated subsheaf of a locally free sheaf $\mathcal S$
is one that is not contained in any larger subsheaf of the same rank and that
it is reflexive. It follows that such a subsheaf,
if it is rank-one, is locally free (see for example, \cite{MP}).
Hence, we can even define $L$ without the dominant condition
on $f$ by setting $m=\dim f(X)$ and replacing $Y$ by a desingularization
of the
algebraic image of $f$ in general. Now $f$ gives rise to an orbifold
rational map in the category of orbifold if an orbifold structure $A$
is imposed on $X$. We denote this orbifold map by $f^\bd$
and the orbifold $X\m A$ by $\xd$ if $A$ is implicit.

\begin{definition} Let $f: X \dra Y$ be a rational map giving rise to an
invertible sheaf $L$ on $X$ as defined above. Let $A$ be an orbifold
structure on $X$ giving rise to an orbifold rational map that we denote by
$$f^\bd= f|_{X\m A}: X\m A \dra Y.$$
Define the vertical part of $A$ with respect to $f$ by
$$
{A\cap f}=\sum \Big\{\ \Big(1-\frac{1}{m(D\cap A)}\Big)D\ \ \Big|\
f(D)\neq f(X),\ D\in \dv'(X)\Big\}.
$$
We set $L_\fd = L(A\cap f)$, which is a $\QQ$-invertible sheaf,  and we
define the
Kodaira dimension of the orbifold rational map $\fd$ by
$$\ka(f, A)=\kappa(\fd):= \kappa(L_\fd).$$
\end{definition}

Recall that a dominant rational map is called almost holomorphic if its
general fibers are
well defined (i.e., do not intersect with the indeterminancy locus).  More
specifically,
the restriction of the second projection to the exceptional locus
of the first projection of
the graph of the map is not dominant. Such a map is
called an almost holomorphic fibration if the general fibers are connected.
Recall also that a fibration is a proper surjective morphism with
connected fibers
while a fiber space is a dominant morphism whose general fibers are
connected.
A dominant rational map is called a rational fibration if it becomes a
fibration
after resolving its indeterminancies.

\begin{definition}
Notation as above, we call the orbifold rational map $f^\bd$ to be 
{\bf\em (base-wise) of general type} (or simply
to be {\bf\em base-general(-typical)}) if $$\ka(f^\bd)=\dim (f)>0,$$
where $\dim (f)$ is given by the dimension of the algebraic image of $f$.
We call
the orbifold $X\m A$ {\bf\em special} if it admits no base-general orbifold
rational map and to be {\bf\em general-typical} or {\bf\em of general type} 
if the identity map restricted to
the orbifold is base general. 
If $f$ is a rational fibration, we say that $f^\bd=f|_{X\m A}$
is {\bf\em base-special} if $X\m A$ has no orbifold rational map that
is base-general and that factors through $f$; We will also consider the
obvious
generalization of this notion to dominant rational maps via Stein
factorization.
If $f$ is an almost holomorphic fibration, we say that $f^\bd$ is {\bf\em
special} (respectively {\bf\em general-typical}) if its
general fiber endowed with the orbifold structure given by the restriction
of $A$
are special orbifolds (repectively orbifolds of general type). 
It should be clear that orbifold structures under
generic restrictions makes sense, see lemma~\ref{lem0}).

\end{definition}
In the case $A$ is reduced and normal crossing,
it is easily seen that these notions are, in the obvious manner,
proper birational invariants of the open subset $X\setminus A$ 
and of the restriction of $f$ to it. Hence, these notions make sense for
quasiprojective varieties and mappings from them and we will so
understand them in this context.\\

This notion of being special corresponds to the same 
``geometric" notion introduced by
Campana in \cite{Ca08} in the case $A_{red}$ is normal crossing and to the
notion given in section 1 in the case $A$ is reduced and normal crossing
by virtue of the following two lemmas respectively, see \cite{Lu02}.
The first of these lemmas is self-evident (with the help of the existence
of diagram~\ref{dm1} in lemma~\ref{B}
as one convenient but not absolutely necessary shortcut).

\begin{lem}\label{log}
Let $L$ be a saturated line subsheaf of $\Omega_X^i$
and $X\m A$ a log-manifold. Then the saturation of $L$
in $\Omega^i(X,\log A)$ is $L(A')$ where $A'$ consists
of components $D$ of $A$ whose normal bundles
$N_D$ over their smooth loci satisfy $N_D^*\wedge L=0$
in  $\Omega_X^{i+1}|_D$. Hence given a dominant map
$\fd: X\m A\dra Y,$  $L_\fd$ is the
saturation of $L_f$ in $\Omega^r(X,\log A)$,
$r=\dim Y$.
\BOX
\end{lem}

\begin{lem}[Bogomolov, Castelnuovo-DeFranchis]
Let $L$ be a saturated line subsheaf of $\Omega^p(X,\log A)$
where $A$ is a normal crossing divisor in $X$. Then
\begin{itemize}
\item[\rm (I)] $\kappa(L)\le p.$
\item[\rm (II)] If $\kappa(L)=p$,
then the Iitaka fibration $I_L$ of $L$ defines
an almost holomorphic fibration
to a projective base $B$ of dimension $p$ and
$I_L^*K_B$ saturates to $L$ in $\Omega^p(X,\log A)$.
In particular, $L=L_{I_L^\bd}$. \BOX
\end{itemize}

\end{lem}

Now let $f: X\ra Y$ be a fibration with $X$ and $Y$ projective
and smooth and let $A$ be an orbifold structure on $X$.
Then the induced orbifold fibration  $f^\bd=f|_{X\m A}$ imposes
an orbifold structure on $Y$ as follows.
Given $D\in \dv'(Y)$, we may write $f^*D=\sum_i m_iD_i$
for $m_i\in \NN$ and
$D_i\in \dv'(X)$. Then we define
the (minimum) multiplicity of $f^\bd$ over $D$ by
$$m(D,f^\bd)=\min \{\ m_i\, m(D_i\cap A)\ |\ f(D_i)=D\ \}.$$

\begin{defn} With the notation as given above,
the $\QQ$-divisor on $Y$
$$D(f^\bd)=D(f|_{X\m A})=D(f,A):=\sum \Big\{\
\Big(1-\frac{1}{m(D,f^\bd)}\Big)D\ \ \Big|\
D\in\dv'(Y)\ \Big\},$$
defines the {\bf\em orbifold base} $Y\m D(\fd)$ of $\fd=f|_{X\m A}$.
\end{defn}
It is immediate that $D(f^\bd)$ is  supported on the union of
$f((A\cap f)_{red})$ with the divisorial part $\Delta(f)$ of the
discriminant locus   of $f$. Note that replacing $f$ by its
composition with a birational morphism $r: \tilde X\ra X$
and imposing the $\infty$ multiplicity along the exceptional
divisor of $r$ while keeping the other orbifold multiplicities
the same does not  change $D(\fd)$. Hence,
although the definition of $D(\fd)$ would no longer
make sense if we allow $f$ to be meromorphic, we
can deal with the problem in a consistent way
(though not always the best way) by
resolving the singularities of $f$ and imposing the $\infty$ multiplicity
along the exceptional divisor of of the resolution. In the case $A$ is
reduced, the same can be achieved by imposing only the $\infty$
multiplicity along the exceptional divisor of $r$ that maps to $A$,
that is, $r^\bd$ gives a proper birational morphism to $X\m A$.
This is always adopted in the case $A$ is reduced.
The following two lemmas
(lemma 3.5 and 3.4 of \cite{Lu02}) are essentially
immediate consequences of the definition.

 \begin{lem}\label{lem0} With the notation as above,
 let $g: Y \ra T$ be a fibration and $h=g\circ f$. Let
 $i: X_t \hra X$, respectively
 $j: Y_t \hra Y$, be the inclusion of the fiber of
 $h$, respectively $g$, above a general point $t\in T.$
 Then $D(\fd)_t:=j^*D(\fd)$ and $A_t:=i^*A$ are
 orbifold structures on the
 nonsingular fibers $Y_t$ and $X_t$ respectively
 and $D(\fd)_t= D(f_t^\bd),$
 where $f_t^\bd=f_t|_{X_t\m A_t}$. That is
 $$ D(f|_{X\m A})|_{Y_t}= D(f_t|_{X_t\m A_t}).$$
Hence $(Y^\bd)_t:=Y_t\m D(\fd)_t$ and $f^\bd_t$ make sense
and $(Y^\bd)_t= Y_t\m D(f^\bd_t)=: (Y_t)^\bd.$

\end{lem}

\noindent
{\bf Proof:} The lemma follows from the fact that $h$,
respectively $g$, and
its restriction to the divisor $R=(A+f^*\Delta(f))_{red}$ in $X$,
respectively the divisorial part of $f(R)_{red}$, are
generically of maximal rank when resticted to their fibers
above $t$ (by Sard's theorem). \BOX\\[2mm]
Hence,
the definition of $D(\fd)$ is well behaved
under generic restrictions.

\begin{lem} \label{comp} Let $\fd$, $g$, $h$ and $A$ be as above, let
$B=D(\fd)=D(f,A)$, $g^\bd= g|_{Y\m B}$ and $\hd=h|_{X\m A}$.
Then $D(\gd)\geq D(h^\bd)$, i.e., $D(g, D(f,A))\geq D(g\circ f, A)$.
If the exceptional part of $A$ with respect to $f$ is reduced or if
$A$ and $B$ are reduced and
$\fd: X\m A\ra Y\m B$ is proper and birational, then
equality holds. \BOX
\end{lem}

The following theorem is
the key fact about special orbifolds used to establish our 
main theorem. It will be used in
the next section.

\begin{theorem}\label{special}
Let $\xd$ be a (smooth) orbifold, $\fd: X^\bd \ra T$ a
special orbifold fibration and $h^\bd: X^\bd \dra Z$ a
base-general orbifold rational map.
Then $h^\bd=k\circ\fd$ for a
rational map $k:T\dra Z$ and $k^\bd:=k|_{T\m D(\fd)}$
is base-general.
In particular, if $T\m D(\fd)$ is special,
then $\fd$ is base special and hence $X^\bd$ is special.
\end{theorem}

This is Proposition~6.5 of \cite{Lu02}, see also \cite{Lu0} and
Chapter 8 of \cite{Ca08}. In view of its importance
here, we reproduce a proof below adapted to our situation.
\\

Recall that a $\QQ$-invertible sheaf is called big if it has maximal Kodaira
dimension.
We first quote two elementary and well-known lemmas concerning the
Kodaira dimension.

\begin{lem}[Kodaira, \cite{Kod}] \label{kod}
Let $H$ and $L$ be invertible sheaves on $X$
with $H$ ample. Then $L$ is big if and only if there is a positive integer
$m$ such that $\dim H^0(L^mH^{-1})\neq 0$.
\end{lem}

\begin{lem}[Easy Addition Law, \cite{Ii}] Let $f : X\ra Y$ be a fibration
with
general fiber $F$ and $L$ an
invertible sheaf on $X$.
Then $$\ka(L)\leq \ka(L|_F)+\dim(Y).$$

\end{lem}

The following is a simplified version of Lemma~5.7 of \cite{Lu02}.

 \begin{lem}\label{gt}
 Consider the following commutative diagram
 of rational maps between complex projective manifolds
\begin{equation*}
\begin{CD}
X @>h>> Z\\
@VV{f \ \ \  \searrow\, g}V @AAwA\\
T@<<u< Y
\end{CD}
\end{equation*}
where $f$ is a fibration, $g$ and $h$ are dominant rational maps
and $u$ and $w$ are morphisms, necessarily surjective.
Let $A$ be an orbifold structure on $X$.
Let $i$ and $j$ be the inclusion
of the general fibers  $X_t:=f^{-1}(t)$ and $Y_t:=u^{-1}(t)$ over $T$.
Let $g_t=g\circ i$ and $h_t=h\circ i$.
Assume that $w\circ j$ is generically finite so that $L_{g_t}=L_{h_t}$.
Then $i^*A$ is an orbifold structure on $X_t$ and
we have with $p=\dim Z=\dim Y$,
$q=\dim Y_t$ that
$$\kappa(h^\bd)- (p-q) \leq \kappa(i^*L_{h^\bd})
\le \kappa(L_{h_t^\bd})=:\kappa(h_t^\bd).$$
In particular, if $h^\bd$ is base-general, then so are $h_t^\bd$ and
$g_t^\bd$ if $\dim Y_t>0$.

\end{lem}

\noindent
{\bf Proof:} It will be clear from our proof that we may assume WLOG,
by taking repeated hyperplane sections of $T$ if necessary, that
$w$ is generically finite.
So we assume this from the start.
Note then that $\dim(T)=p-q$ so that the first inequality above
follows from the easy addition law of Kodaira dimension.
To obtain the second inequality and thus the lemma, the
easily verified fact that
$$i^*(A\cap h)=(i^*A)\cap (h\circ i) = (i^*A)\cap h_t =(i^*A)\cap g_t$$
allows us to deduce it from an inclusion of $i^*L_h$ in $L_{h_t}=L_{g_t}$
that can be seen as follows. The conormal short exact sequence
on $X_t$
	$$0\ra N^*_{X_t}\ra \Omega_X|_{X_t} \ra \Omega_{X_t} \ra 0$$
gives rise to a natural sheaf morphism from $i^*\Omega_X^p=\Omega_X^p|_{X_t}$
to the factor $\Omega^q_{X_t}\otimes \Lambda^{p-q} N^*_{X_t}$ in its quotient
filtration.
Now, over the Zariski open set $U$ of $X_t$ where $g_t$ and $h_t$ are
defined,
$i^*g^*=g_t^*$ gives a map from the same short exact sequence on $Y_t$ to
that
of $X_t$. So it does the same for the corresponding
natural sheaf morphism on $Y_t$  to that of
$X_t$. Thus we obtain a commutative diagram over $U$:
\begin{equation*}
\begin{CD}
i^*h^*K_Z@>>>\Omega_X^p|_{X_t} @>>>
\Omega_{X_t}^q\otimes\Lambda^{p-q}N_{X_t}^*\\
@VV\delta V@AAA@AAA\\
i^*g^*K_Y@>\sim>>g_t^*(K_Y|_{Y_y})@>>> g_t^*(K_{Y_t}\otimes \det N^*_{Y_t})
\end{CD}
\end{equation*}
where $\delta$ is induced from the inclusion $w^*K_Z\hra K_Y$.
As both $\det N^*{Y_t}$ and $\Lambda^{p-q}N_{X_t}^*=\det N_{X_t}^*$ are
trivial invertible sheaves by construction, we see that $g_t^*K_{Y_t}$ has
the
same image in $\Omega^q_{X_t}$ as that of $i^*h^*K_Z$ over a Zariski open
subset of $X_t$. As the former saturates
to $L_{g_t}$ in $\Omega^q_{X_t}$, we see that $i^*L_h\hra L_{g_t}$ as
required.
\BOX\\

\noindent {\bf Proof of Theorem~\ref{special}}: Let $ g_0=(f,h): X\dra
T\times Z$
and $Y_0$ its image. Let $r: Y\ra Y_0$ be a resolution of singularities of
$Y_0$
and let $g=g_0\circ r^{-1}$, which is a rational map in general. Let $u$
and $w$
be r composed with the projections of $T\times Z$ to $T$ and $Z$
respectively.
Then we are in situation of Lemma~\ref{gt} with $h^\bd$ base-general. We
first note that the general fibers of $u$ are connected by construction,
being
images of the fibers of the special fibration $f$.
As the general fibers of $f$ are special, our lemma implies that $Y_t$ for
the general $t\in T$ are points. It follows that $u$ is birational so that
$Y_0$
form the graph of a rational map $k: T\dra Z$ and $h=k\circ f$. The theorem
now follows directly from the following elementary lemma,
which is a simplification of
Proposition~3.19 of \cite{Lu02}. \BOX

 \begin{lem} \label{lem}
 Let $f : X\ra T$ be a fibration,  $k: T\dra Z$ a rational map
 and $h=k\circ f$.
 Let $A$ be an orbifold structure on $X$ inducing the
 orbifold maps $\fd=f|_{X\m A}$ and $h^\bd=h|_{X\m A}$.
 Let $B=D(\fd)$ be the orbifold structure
 on $T$ imposed by $\fd$ and let $k^\bd=k|_{T\m B}$ be the induced orbifold
 map from $T$. Then $$\ka(h^\bd)\leq \ka(k^\bd).$$
In particular, if $h^\bd$ is base-general, then so is $k^\bd$.
\end{lem}

\noindent
{\bf Proof:} We only prove the lemma in the case we are using in the paper
where $A$ and $B$ are reduced normal crossing divisors;
more specifically, in the case $A$ is reduced normal crossing,
$\fd$ is a special fibration and and $B=D(\fd)$ is just the standard boundary
divisor of the compactification of a semi-Abelian variety.
In this case,
both $X\m A$ and $Y\m B$ are log-manifold and so $L_{\fd}$ and $L_{k^\bd}$
can be considered respectively as invertible subsheaves of
$\Omega(X,\log A)$ and $\Omega(T,\log B)$
by Lemma~\ref{log}. We note that, outside the exceptional divisor $E(f)$
of $f$, $f$ is a log-morphism (i.e., $f^{-1}(B)\subset A$) and so
gives an inclusion of sheaves $f^*\Omega(T,\log B)\hra \Omega(X,\log A)$
there and it is actually a vector bundle inclusion on a Zariski open subset
of $f^{-1}(B)$ (\cite{Ii}).
Thus, we have an inclusion $$f^*L_{k^\bd} \hra L_{h^\bd}$$ outside $E(f)$
that is an equality  on a Zariski open subset
of $f^{-1}(B)$. This equality extends to the open subset outside
$A\  (\supset f^{-1}(B))$
where $f$ is smooth. But by our definition of the multiplicity that gives the
orbifold base, this open subset surjects to the complement of a codimension
two subet of $T\m B$. Hence
$H^0(L_{h^\bd}^l) \hra H^0(f^*L_{k^\bd}^l)=H^0(L_{k^\bd}^l)$ for all
positive integer $l$ by the Hartog extension theorem. \BOX\\[2mm]
We remark that in our case at hand, $Y\m B$, being a semi-Abelian
variety with an equivariant compactification $Y$ (see the next section
for the definition and basic facts),
has trivial log-cotangent sheaf. Hence $\ka(L_{k^\bd})\leq 0$ and $Y\m B$
is a special orbifold. \\

We now address the very important question of when is the base Kodaira
dimension of an orbifold fibration equal to the
Kodaira dimension of its orbifold base. 
The question was posed by Campana in \cite{Ca02} for which 
he gave a partial answer in the case the base has positive 
Kodaira dimension. We have also given a partial answer
in lemma~2.2 and 2.4 of \cite{Lu02} 
which showed at the same time the equivalence of our approach
to that of Campana's. It is this latter that we give below but restricted
here for simplicity to the context of log-manifold.

\begin{lemma}\label{B}
Let $f: X\ra Y$ be a fibration where $X$ and $Y$ are
complex projective manifolds, $A$ a normal crossing
divisor on $X$ and $\fd=f|_{X\m A}$. Then
$\kappa(\fd) \leq \kappa(Y\m D(\fd)).$
Also, one can find a commutative
diagram of morphims between complex projective manifolds
\begin{equation}\label{dm1}
\begin{CD}
X' @>v>> X\\
@Vf'VV @VVfV\\
Y' @>>u> Y
\end{CD}
\end{equation}
with $u,v$ birational and onto such that
$E(f'\circ v^{-1})=\emptyset$ and 
$A'=v^{-1}(A)$ is normal crossing.
Let  $f'^\bd=f'|_{X'\m A'}$.
Then $v$ induces a
proper birational morphism
$X\setminus A\ra X'\setminus A'$ and
$\ka(\fd)=\ka(f'^\bd)$.
If $m$ is divisible by the multiplicities of $D(f'^\bd)$, then
$$H^0(X, L_{f'^\bd}^m)= H^0( Y', K_{Y'\m D(f'^\bd)}^m)\ \ \text{and }\ \
\kappa(f'^\bd)= \kappa(Y'\m D(f'^\bd)).$$
\end{lemma}

\noindent
{\bf Proof:} The first statement
follows from lemma~\ref{lem} by letting $k$
be the identity map there.

The construction of a birationally
equivalent fibration as given by the commutative
diagram with the above property is
achieved by resolving the singularities of the flattening of
$f$, which exists by \cite{R, GR}, and in such way that the inverse image of 
$A$ is normal crossing, which is possible by \cite{Hi}.
As $u$ is birational, $v^*L_\fd \hra L_{f'|_{X'\m A'}}$
by lemma~\ref{log} and hence
$\ka( L_\fd)\leq \ka(L_{f'|_{X'\m A'}})$. The reverse inequality
follows  from  lemma~\ref{lem}.

For the last statement, we have
(with $r=\dim Y'$) as before the inclusion
\begin{equation}\label{sat1}
{f'}^*K_{Y'}(D(f'^\bd))^m\hra L_{f'^\bd}^m\
\Big(\hra (\Omega^r_{X'}(\log A'))^{\otimes m}\Big)
\end{equation}
outside $O\cup E(f')$ where $O$ is a subset of
$X'$ of codimension two or higher
contained above the discriminant $\Delta(f')$ of $f'$ and
$m$ is a positive multiple of all relevant multiplicities.
Moreover, this inclusion is an isomorphism on an open subset
of $X'$ that surjects to the complement of a
subset of $Y'$ of codimension two or higher. Hence
$H^0(Y', K_{Y'}(D(f'^\bd))^m) \hra H^0(X, L_\fd^m)
=H^0(X', L_{f'^\bd}^m)$
by the Hartog extension theorem applied to $X$ and the
reverse inclusion by the Hartog extension theorem applied
to $Y'$. $\BOX$\\

We give below generalizations 
to the relative setting of lemma~\ref{B} and theorem~\ref{special}.
They are used to extend our main theorem
but are otherwise not needed for its proof.

 \begin{lem}\label{lem1}
 With the setup as in lemma~\ref{B} and with all elements of the
 commutative diagram (\ref{dm1}) as given there,
 let $g':Y'\dra Z'$ be a dominant rational map.
 Then $$\ka(g', D(f'|_{X'\m A'}))=\ka(g'\circ f',A').$$
\end{lem}
\noindent
{\bf Proof:} The proof is the same as that of lemma~\ref{B},
replacing $K_{Y'}(D(f'^\bd))$ by $L_{g'|_{Y'\m D(f'^\bd)}}$.
\BOX\\

\begin{prop}\label{s/s}
Let the setup be as in lemma~\ref{lem0}, i.e.,
$f:X\ra Y$ and $g:Y\ra T$ are fibrations with orbifold structures
$A$ on $X$ and $D(f, A)$ on $Y$ and $t$ a general point on $T$.
Assume $A$ is a normal crossing divisor.
If $\fd$ and $Y_t^\bd=Y_t\m D(f|_{X^\bd_t})$ are special,
then so is $X^\bd_t:=X_t\m A_t$. If $T\m D(g\circ f, A)$ is special,
then $\fd$ is base-special if so is $f_t^\bd$.
\end{prop}
{\bf Proof:} The last statement follows by noting that $\fd$ is 
base-special if so is $f'^\bd$, in the notation of lemma~\ref{B}, and
$f'^\bd$ is base-special by lemma~\ref{lem1} and so 
theorem~\ref{special} applies. \BOX

\section{Structure of the quasi-Albanese map}

Let $X$ be a complex quasi-projective manifold, 
$T$ a semi-Abelian variety and
$u: X\ra T$ an algebraic morphism. Let $\bar T$ be
a smooth equivariant compactification of $T$, i.e., $\bar T$ is smooth
and admits
an algebraic action by $T$ -- an example being
the compactification of $T$ in the exact sequence~(\ref{2})
via the compactification $\CC^k\subset (\PP^1)^k$.
Then one can observe, see \cite{NW},
that $\bar T\setminus T$
is a normal crossing divisor and that $\bar\Omega_T$
is a trivial bundle over $\bar T$ (via simultaneous equivariant resolution
of singularities for example).
By the resolution of singularity theorem (see \cite{BM, Hi}),
there is a compactification $\bar X$
of $X$ with normal crossing boundary divisor $A$
such that  $u$ extends to a morphism
$\bar u: \bar X \ra \bar T$.

\begin{defn} We call $\bar u$ as above a {\bf\em natural compactification} of $u$.
We will set $\bar u^\bd=\bar u|_{\bar X\m A}$
and, in the case $u$ is a fiber space (i.e., $\bar u$ is a fibration), we set
$D(u) := D(\bar u^\bd)|_{X}\ .$
\end{defn}
We note that $D(u)$ is a $\QQ$-divisor on $X$ that is independent of the
natural compactification $\bar u$ of $u$ chosen since
two such
compactifications are always dominated by 
a third such compactification. By the same token, the notions of being
special, being general-typical, being base-special and
being base-general(-typical) are well-defined for $u$ (independent of
the natural compactifications).

\begin{defn}[\cite{NWY2}] Let $D$ be an algebraic subset or a 
$\QQ$-divisor in $T$. We define
${\rm St}(D)$ to be the identity component of $\{a\in T: a+D=D\}$
which is easily verified to be a subgroup, even a semi-Abelian
subvariety. Given a compactification $\bar T$ of $T$, we define
$\bar D$ to be the Zariski closure of $D$ in $\bar T$. 
\end{defn}

The purpose of this section is to establish the following proposition
but restricted to the situation
we are in.

\begin{prop}\label{str}

Let $X$ be a complex quasi-projective manifold with (quasi-)Albanese map
$ f:  X \ra T_0={\rm Alb}(X)$. If $f$ is not base special,
then there is a proper semi-abelian subvariety $T'$ of $T_0$
such that if $g: T_0\ra T=T_0/T'$
is the quotient map, then the orbifold base of 
$h=g\circ f$ is of general type and of positive dimension; More specifically, 
if $\bar e: \bar Z\ra \bar T$ is the finite map factor in 
the Stein factorization of $\bar h$, then $\dim \bar Z>0$ and either 
$Z=\bar e^{-1}(T)$ is of general type or $h$ is a fiber space
(i.e., $\bar h$ a fibration)
and $K_{\bar T}(D(\bar h, A))={\mathcal O}(\bar D)$ 
is big for one and hence for
all natural compactification  $\bar h: \bar X\ra \bar T$ of $h$
with $A=\bar X\setminus X$ and $D=D(h)$.    \\

Suppose that $f$ is a fiber space and let $h_0=f$. 
Consider the following inductive definition. With
$h_{i-1}$ and $T_{i-1}$ defined, we 
define $D_{i-1}=D(h_{i-1})$,
$T''_i={\rm St}(D_{i-1})$, $T_i=T_{i-1}/T''_i$  and
$h_i=\ga_i\circ h_{i-1}$ with $\ga_i:T_{i-1}\ra T_i$.
Then this process terminates at the $l$-th stage
for $l$ such that ${\rm St}(D_l)=\{0\}$. Setting $h=h_l$
$D=D(h)_{red}$ and $T=T_l$, 
we find two possibilities:

\begin{itemize}

\item[(i)] $\dim T>0$ and 
${\mathcal O}(\bar D)$ on $\bar T$ is big 
for any equivariant compactification $\bar T$ of $T$.

\item[(ii)] $T$ reduces to a point and $f$ is base-special.\\

\end{itemize}

\end{prop}
In order to prove this proposition, we
first recall the structure theorem of Kawamata \cite{Ka},
Kawamata-Viehweg \cite{ KV} and
Ueno \cite{Ue} concerning the Albanese map:

\begin{prop}\label{KV}
Let $Z$ be a normal quasiprojective variety with a 
finite morphism to a semi-Abelian variety $T$.
Then there is a finite extension ${\rm St}(Z)$ of 
an abelian subvariety of $T$ whose natural action 
on $T$ lifts to $Z$ and 
$Z/{\rm St}(Z)$ is of general type of dimension
$\bar \ka(Z)=\dim Z- \dim {\rm St}(Z)\geq 0$. In particular,
$f$ is an \'etale covering map over a translate of
a semi-Abelian subvariety of $T$ if and only if $\bar\ka(Z)=0$
if and only if $Z$ is semi-Abelian. \BOX
\end{prop}

Since the image of the Albanese map  $f=\alpha_X$
of $X$ generates Alb$(X)$ (whose fact is equivalent to
the universal property of the Albanese map),
applying the above proposition
to the finite map factor of the Stein factorization of $f$
yields directly proposition~\ref{str} in
the case $f$ is not a fiber space.
In the case $f$ is a fiber space, the first
paragraph of the proposition follows directly from the other 
part of the proposition. 
Hence, it remains to establish the two cases  
(i) and (ii) of the proposition to end the
proof of the proposition.
We need the following lemma, well-known
in the case of Abelian varieties but generalized to the semi-Abelian
case by Proposition 3.9 of \cite{NWY2}.

\begin{lem}\label{kabig}
 Let $D$ be an effective $\QQ$-divisor in
 a semi-Abelian variety $T$,
 $\bar T$ an equivariant compactification of $T$ and
 $\bar D$ the Zariski closure of $D$ in $\bar T$. 
 The the conditions ${\rm St}(D)=\{0\}$, 
 ${\mathcal O}(\bar D)$ being big, and ${\mathcal O}(\bar D_{red})$ 
 being big are equivalent.
\end{lem}

We remark that for an effective $\QQ$-divisor $D$, 
there exists $m>0$ such that
$\frac{1}{m} D_{red}\leq D\leq mD_{red}$
and St$(D)={\rm St}(D_{red})$.
Hence, the lemma follows trivially from the weaker
assumption that $D$ is a an effective 
(even a reduced) divisor.
We remark also that one can give a direct proof of the above lemma
using the original
arguments of the above structure theorem, or by applying
the above theorem to suitable ramified covers of $T$. \\

\noindent
{\bf Proof of proposition~\ref{str}:} 
The claim for case (i) is just lemma~\ref{kabig}. 
Let $g_i=\ga_i\circ\dots\circ\ga_1$.  To prove (ii), it
suffices to prove for general $t_i\in T_i$
that  $f_{t_i}=f|_{h_i^{-1}(t_i)}$ is base special for all $i$. But this follows
by induction as follows. It is clear for $i=0$. Assume that it is true for
$f_{t_{i-1}}$.
The orbifold base of $h_{i-1}|_{t_i}=g_{i-1}\circ f_{t_i}$, 
being the restriction of that of $h_{i-1}$ to $\ga_i^{-1}(t_i)$,
is simply $\ga_i^{-1}(t_i)$ by lemma~\ref{lem0}, which is a 
general translate
of the semi-Abelian subvariety $T''_i$ and is thus special.
As
$f_{t_i}|_{t_{i-1}}=f_{t_{i-1}}$ is base-special by the induction 
hypothesis, lemma~\ref{s/s} shows that $f_{t_i}$ is base-special
as required. \BOX

\section{Implication of Zariski dense entire curves}

We first recall some relevant definitions and facts from Nevanlinna theory.
We will follow section 2 of \cite{NWY2}. Let $T$ be a complex
manifold, $\omega$ a real smooth $(1,1)$-form and $\ga : \CC\ra T$
a holomorphic map. Then the order function of $\ga$ with respect
to $\omega$ is defined by
\begin{equation}
\label{2.1}
T_\ga(r; \omega)=\int_1^r\frac{dt}{t}\int_{|z|<t} f^*\omega
\qquad (r>1). 
\end{equation}
If $T$ is K\"ahler and $\omega$, $\omega'$ are $d$-closed
real (1,1)-forms in the same cohomology class $[\omega]$, then
$$
T_\ga(r; \omega)=T_\ga(r; \omega')+O(1).
$$
Hence we may set, up to $O(1)$-terms,
\begin{equation}
\label{2.2}
T_\ga(r; [\omega])=T_\ga(r; \omega).
\end{equation}
Let $L \to T$ be a hermitian line bundle. 
As its Chern class is a real $(1,1)$-class,
we may set
$$
T_\ga(r; L)=T_\ga(r; c_1(L)).
$$
We will denote by ${\mathcal O}_T(D)$ the line bundle
determined by a divisor $D$ on $T$ via a standard abuse of notation
and set $T_\ga(r; D)=T_\ga(r; {\mathcal O}_T(D))$. 
By using a standard Weil function for
to a subscheme $W$ of $T$ (see section 2 of \cite{NWY2}), 
we can define $T_\ga(r; W)$ similarly, which we will also denote
by $T_\ga(r; {\mathcal I}_W)$ where ${\mathcal I}_W$ is the ideal subsheaf of 
${\mathcal O}_T$ defining $W$. 
This is because, ideal sheaves pulls back to the same
on $\CC$ which are then ideal sheaves defining effective divisors on
$\CC$. Henceforth, we will identify ideal sheaves with their 
effective divisors on $\CC$.\\

\noindent
Let $E=\sum_{z\in \CC} ({\rm ord}_z E) z$ 
be an effective divisor on $\CC$,
$S\subset \CC$ and $l\in \NN\cup \{\infty\}$. 
Then the above sum is a sum over a discrete subset of $\CC$.
Hence the sum is finite 
when restricted to the disk $\DD_t$ of radius $t>0$
and so $$n(t;E)=\deg_{\DD_t} E:= \sum_{z\in \DD_t} {\rm ord}_z E$$
makes sense. We define the restriction of $E$ to $S$ truncated to
order $l$ by $$E_{S,l}=\sum_{z\in S} \min({\rm ord}_z E, l) z$$
and set $n_l(t;E,S)=n(t; E_{S,l})$.
Then the counting functions of $E$ with, respectively without,
truncation to order $l$ are given by
$$
N_l(r; E, S) =\int_1^r \frac{n_l(t; E, S)}{t}dt\ ,\ \ \ N_l(r;E)=N_l(r;E,\CC),
$$
respectively $N(r; E)=N_\infty(r; E)$.\\

A well-known consequence of these definitions via the classical
Jensen formula is the First Main Theorem:
\begin{equation}\label{m}
N(r;\ga^*{\mathcal I}_W)\leq T_\ga(r; W) + O(1),
\end{equation}
where $W$ is a subscheme of $T$ 
and $\ga^*{\mathcal I}_W$, as an ideal subsheaf of ${\mathcal O}_\CC$,
is identified with a divisor on $\CC$. By the linearity of 
$n_l$ and hence of $N_l$ with respect to the third variable, 
if $T$ is a disjoint union of $U$ and $V$,
then
$$N_l(r; \ga^*{\mathcal I}_W)=N_l(r; \ga^*{\mathcal I}_W,\ga^{-1}(U))
+N_l(r; \ga^*{\mathcal I}, \ga^{-1}(V)).$$ 
Also if $Z\subset E$ is an inclusion of reduced algebraic
subsets of $T$ and $Z=\ $supp$\,\mathcal I$, then
\begin{equation}\label{mm}
N(r; \ga^*{\mathcal I}_Z)\geq N(r; (\ga^*{\mathcal I}_Z)_{red})
=N_1(r; \ga^*{\mathcal I}_Z)= N_1(r; \ga^*{\mathcal I}_E,\ga^{-1}(Z)).
\end{equation}
By definition, we also have for
$D$ a strictly effective divisor on $T$
that $T_\ga(r; D) \geq 0$ if the image
of $\ga$ is not in $D_{red}$.
This fact along with Kodaira's lemma~(\ref{kod})
and linearity of $T_\ga$ with respect to to the second variable 
yields easily the following
(which is lemma~2.3 of NWY2):

\begin{lemma}\label{big}
Suppose that $T$ is a complex projective manifold with a
big divisor $A$. Then $$T_\ga(r)= O(T_\ga(r; A)).$$
\end{lemma}
Here by convention $T_\ga(r):= T_\ga(r; H)$ for an ample divisor $H$ on $T$.\\

We will need the following two theorems of Noguchi-Winkelmann-Yamanoi. 

\begin{theorem}[\cite{NWY3}]\label{NYWIII}
Let $X$ be a normal complex quasi-projective variety admitting
a Zariski dense entire holomorphic curve.
Let $f: X\ra T$ be a finite morphism to a semi-Abelian variety.
Then $f$ is an \'etale covering morphism.
\end{theorem}

\begin{theorem}[\cite{NWY2}]\label{NWYII}
Let $T$ be a semi-abelian variety and $\ga:\CC\ra T$ a holomorphic
map with Zariski dense image. Let $E$ be a divisor on $T$ and $\bar E$
be its Zariski closure in a equivariant compactification of $T$. Let $\mathcal I$ be
an ideal subsheaf of ${\mathcal O}_T$ such that ${\mathcal O}_T/{\mathcal I}$ is
supported on a codimension-two subvariety of $T$.
Then we have :
\begin{equation} \label{mt10}
N(r;\ga^*{\mathcal I}) \leq \epsilon T_\ga(r)||_{\epsilon} \: {\rm for\: all}\:
\epsilon >0
\end{equation}
\begin{equation} \label{mt20}
T_\ga(r; \bar E) \leq N_1(r; \ga^*E) +
\epsilon T_\ga(r; \bar E)||_{\epsilon} \: {\rm for\: all}\: \epsilon >0
\end{equation}

\end{theorem}
Here ``$||_\epsilon$'' stands for the inequality to hold
for every  $r>1$ outside a Borel set of finite Lebesgue measure
that depend on $\epsilon$.\\

\begin{cor}\label{N3}
Let $X$ be a smooth complex quasi-projective variety admitting
a Zariski dense entire holomorphic curve.
Then the (quasi-)Albanese map of $X$ is a fiber space.
\end{cor}

\noindent
{\bf Proof of corollary:}  This is an easy argument applying  
theorem~\ref{NYWIII} to the
Stein factorization of a natural
compactification of $f$. \BOX\\

Our main theorem in this paper is a direct consequence of the
following proposition.

\begin{prop}\label{last}
Let $X$ be a smooth complex quasi-projective variety admitting
a Zariski dense entire holomorphic curve.
Then its (quasi-)Albanese map $\alpha: X\ra T_0={\rm Alb}(X)$ 
is a fiber space and
$\alpha$ is base special.
\end{prop}

\noindent
{\bf Proof:} The first part of the proposition is just corollary~\ref{N3}. 
We may thus assume that $\alpha$ is a fiber space but not base special.
The proposition is proved once we reach a contradiction with
the existence of an entire map $\ga_0: \CC\ra X$ with Zariski dense
image in $T_0$.
By proposition~\ref{str}(i), there is a quotient morphism $g: T_0\ra T$
of semi-Abelian varieties such that the
fiber space $h=g\circ f: X\ra T$ induces an orbifold base of general
type with $\dim T>0$, that is, 
given any equivariant compactification
$\bar T$ of $T$, ${\mathcal O}(\bar D)$ is big where $D=D(h)_{red}$.
Let $\ga=h\circ \ga_0$.
As $\ga$ has Zariski dense image in $T$, theorem~\ref{NWYII} gives
\begin{equation} \label{mt2}
T_\ga(r;\bar D) \leq N_1(r; \ga^*D) +
\epsilon T_\ga(r; \bar D)||_{\epsilon} \: {\rm for\: all}\: \epsilon >0
\end{equation}
From the definition of $D(h)$, we see that outside the
exceptional locus $E$ of $h$, the effective divisor
$D_0=h^*(D)$ is nowhere-reduced
along $h^{-1}(D)$. Since $Z=h(E)$ is of co-dimension two or higher,
we have from theorem~\ref{NWYII} that
\begin{equation} \label{mt1}
N(r;\ga^*Z) \leq \epsilon T_\ga(r)||_{\epsilon} \: {\rm for\: all}\:
\epsilon >0.
\end{equation}
As $\ga^*(D)=\ga_0^*(D_0)$ is then an effective divisor that
is nowhere-reduced outside $\ga^{-1}(Z)$, we have by (\ref{mm}) that
\begin{equation} \label{inter}
N(r; \ga^*(D)) \geq 2 N_1(r; \ga^*D, \ga^{-1}(T\setminus Z))
\geq 2 N_1(r; \ga^*(D)) - 2 N(r; \ga^*{\mathcal I}_Z).
\end{equation}
Also, since $\bar D$ is big,
we have by lemma~\ref{big},
that $T_\ga(r) = O(T_\ga(r; \bar D)$. That is,
there exists a constant $C>0$
(depending on the ample divisor used to define $T_\ga(r)$) such that
$$T_\ga(r) \leq C T_\ga(r; \bar D).$$
Couple all these with the First Main Theorem, equation (\ref{m}), gives
$$T_\ga(r; \bar D) +O(1)\geq N(r; \ga^*D)
\geq 2 N_1(r; \ga^*D) - 2 N(r; \ga^*Z) $$ $$\geq
(1 - \epsilon) 2T_\ga(r; \bar D)
-2\epsilon C T_\ga(r; \bar D)||_{\epsilon}$$
$$\geq 2(1- \epsilon - \epsilon C) T_\ga(r; \bar D)||_{\epsilon},$$
valid for all $\epsilon>0$.
This gives a contradiction for $\epsilon >0$ sufficiently small. \BOX

\end{document}